\documentclass[12pt,a4]{article}

\usepackage{graphicx}
\usepackage{epsfig}
\usepackage{psfrag}
\usepackage{amsfonts}
\usepackage{amsmath}
\usepackage{amssymb}
\usepackage{color}
\usepackage{bbm}

\newtheorem{lemma}{Lemma}[section]
\newtheorem{theorem}{Theorem}[section]

\newtheorem{assumption}{Assumption}[section]
\newtheorem{algorithm}{\underline{Algorithm}}[section]

\newlength{\subtextwidth}
\newcommand{\diag}{\operatorname{diag}}
\newcommand{\blkdiag}{\operatorname{blkdiag}}
\def\proof{\noindent{\bf{\em Proof:}\ \ }}
\def\QED{\mbox{\rule[0pt]{1.5ex}{1.5ex}}}
\def\endproof{\hspace*{\fill}~\QED\par\endtrivlist\unskip}
\normalbaselineskip=\baselineskip
\newenvironment{eqtarraya}{\begin {displaymath}%
\begin {array}{rcl@{\extracolsep{1cm}}l}}{\end{array}\end{displaymath}}
\newenvironment{eqtarrayb}{\begin {displaymath}%
\begin {array}{r@{\extracolsep{.1cm}}rcl}}{\end{array}\end{displaymath}}
\newenvironment{eqtarrayc}{\begin {displaymath}%
\begin {array}{llll@{\extracolsep{.5cm}}l}}{\end{array}\end{displaymath}}
\def\dsty{\displaystyle}

\def\v0{\hbox{\bf 0}}

\def\R{\hbox{\bf R}}

\def\quad{\hspace{1 cm}}
\def\qquad{\hspace{.5 cm}}

\def\day{\partial}
\def\diag{\mbox{diag}}

\def\bn{\begin{enumerate}}
\def\en{\end{enumerate}}
\def\bq{\begin{eqnarray}}
\def\eq{\end{eqnarray}}
\def\bqn{\begin{eqnarray}}
\def\eqn{\end{eqnarray}}
\def\bqta{\begin{eqtarraya}}
\def\eqta{\end{eqtarraya}}
\def\bqtb{\begin{eqtarrayb}}
\def\eqtb{\end{eqtarrayb}}
\def\bqtc{\begin{eqtarrayc}}
\def\eqtc{\end{eqtarrayc}}
\def\be{\begin{equation}}
\def\ee{\end{equation}}
\def\bea{\begin{eqnarray}}
\def\eea{\end{eqnarray}}
\def\beann{\begin{eqnarray*}}
\def\eeann{\end{eqnarray*}}
\def\bsea{\begin{subeqnarray}}
\def\esea{\end{subeqnarray}}
\def\bmat{\left[ \begin{array}}
\def\emat{\end{array} \right]}
\def\nn{\nonumber}
\def\bd{\begin{displaymath}}
\def\ed{\end{displaymath}}
\def\bi{\begin{itemize}}
\def\ei{\end{itemize}}
\def\bem{\begin{em}}
\def\eem{\end{em}}
\newcommand{\real}{{\mathbb{R}}}
\def\dots{\hbox{ ...}}

\begin{document}

\title {Nonovershooting output regulation  for feedback linearisable  MIMO systems}
\author{Robert Schmid  \thanks{Robert Schmid is  with  the  Department of Electrical and Electronic Engineering, University of Melbourne.  email: rschmid@unimelb.edu.au } }
\date{ \ }
\maketitle

\vspace{-1.5cm}
\begin{abstract}  We  consider the  problem of   designing  a state feedback  control law   to achieve nonovershooting  tracking for  feedback linearisable multiple-input multiple-output nonlinear systems. The reference  signal is assumed to be  obtained  from a  linear exosystem.  The design method adopts  known methods for  the  nonovershooting  tracking  of linear systems  and   extends them to the output  regulation of   feedback  linearisable nonlinear systems.  The method accommodates arbitrary initial  conditions and  offers design choice of the tracking convergence speed.
\end{abstract}

\section{Introduction}
 The problem of  ensuring the system response tracks a desired reference signal  with zero steady-state error and desirable transient  performance  is one  of the classic problems  of  control  systems  theory.   The  twin performance of  objectives  of obtaining a rapid response  while avoiding  or  minimising overshoot  have traditionally  been viewed as competing  objectives,  with control methods seeking a suitable trade-off between the two \cite{Chen-03}.   The problem of entirely avoiding  overshoot has most been considered for linear  systems.  The early works
  \cite{Lin-97}-\cite{Bement-04} considered the step  response of  linear  single-input single output (SISO) systems.  The  problem of avoiding  overshoot  in the step response of a  multiple-input multiple  output (MIMO) linear system was first considered in \cite{SN10},  and this was extended  to  the nonovershooting tracking  of  time-varying   signals   in \cite{SN14}.   Recently \cite{Ntog-16} gave  necessary and sufficient conditions for the  existence  of  static state feedback  controller to  achieve a  monotonic step response,  for  MIMO linear systems.

To  date only a little consideration  has  been  given to  the  problem of nonovershoot for nonlinear systems,  the  principal contributions  being \cite{Krstic-06} and  \cite{Zhu-13}.  Both of these papers  pointed to  notable  `gaps' in the  linear control systems literature.
The  authors  of  \cite{Krstic-06}  contemplated approaching the problem of  designing nonovershooting controllers for  feedback linearisable systems \cite{Shastry-99}   by  first  converting the nonlinear system  into a linear system in chain of integrators form  via a suitable coordinate change and feedback, and then using ''some standard linear technique'' to  achieve a  nonovershooting response. However the authors noted that the coordinate change would, in general,  yield a non-zero initial condition, and (at their time of writing) the techniques available for linear systems assumed zero initial conditions. Thus  the feedback linearization approach was  not   adopted in \cite{Krstic-06}.  Instead the authors  considered
  nonlinear systems in  strict feedback form with the first state as output, and employed a modified back-stepping technique  to obtain a closed-loop  system with a cascaded structure  that was shown to  have   a nonovershooting response, if the feedback  gains were chosen appropriately.

Continuing this  line of research, the authors of \cite{Zhu-13} noted that converting the nonlinear system into chain-of-integrator normal form  made the  problem of  avoiding  overshoot relatively more challenging than for other conventional linear systems,  as integrators are known to  be an important source of overshoot. To address this problem, the authors combined the feedback linearisation method with a global coordinate transformation that they applied to the chain-of-integrators system. The problem  of  ensuring the error signal  did not change sign
 (and hence yield nonovershooting tracking) was shown to be equivalent to ensuring a certain  closed-loop transfer function had a non-positive  impulse response. For chain-of-integrator systems with relative degree of at most four, the authors established sets of feedback  gains  to  ensure a nonovershooting response from an arbitrary initial condition.

 The present work advances upon these two papers in several directions. Following   the  direction  contemplated in \cite{Krstic-06}  and utilised in \cite{Zhu-13},  we assume the system  may  be rendered in chain-of-integrators form by feedback linearisation;  the resulting   zero  dynamics are also  required to be stable. However, unlike  \cite{Zhu-13},   we do  not employ  an additional coordinate change to the chain-of-integrators system. Instead we adapt the nonovershooting control  methods of \cite{SN10} and  \cite{SN14} to linear chain-of-integrator systems.     These methods involve selecting candidate sets  of closed-loop eigenvalues and applying simple analytic tests to  determine the shape of the closed-loop system response arising from the specified initial condition.  When the tests  indicate a suitable set of closed-loop eigenvalues has been found,  the  pole-placing feedback matrix  that delivers the nonovershooting response can  easily be computed from Moore's method \cite{Moore76}.  While \cite{Krstic-06} and \cite{Zhu-13} considered only  SISO systems,  the methods  of \cite{SN10} and  \cite{SN14} are  presented for  MIMO systems and  hence our  presentation will  also be given in a MIMO setting.

   The paper  is  organised as follows. In Section 2 we introduce the problem of  nonovershooting  output regulation for  nonlinear systems,  and  present the  system assumptions required for feedback linearisation. Section 3 revisits some results  from  \cite{SN10} and  \cite{SN14}  and adapts them to  linear systems in chain-of-integrators form. Section 4 presents the main results of the paper on the controller design method. Section 5  considers an  example system that was also considered in \cite{Krstic-06} and \cite{Zhu-13}, to enable comparisons of the design methods and their relative performance.

\section{Output regulation of feedback linearizable MIMO nonlinear systems}

We consider  an affine nonlinear square MIMO  system $\Sigma_{nonlin}$  in the form
\bea
\Sigma_{nonlin}:\; \left\{ \begin{array}{rcl}
\dot x    & = &  f(x) + g(x) u(x), \quad x_0 = x(0)  \\
 y       & = &  h(x)
  \end{array} \right.   \label{nonlineq}
\eea
where  $x  \in \real^n$, $u  \in \real^p$,   $y \in \R^p$, and  $f$,   $g$   and  $h$ are smooth  vector fields.
  The problem of {\it nonovershooting  output regulation}  is to  find a state feedback  control law $u(x)$ that stabilizes the closed-loop  system and ensures the system output tracks a reference signal $r \in \real^p$ without  overshoot; thus $e(t) = r(t) - y(t) \rightarrow 0$ without changing sign in  all  components.  We  shall  assume the reference signal is obtained as the output of a linear  exosystem
      \be
    \Sigma_{exo}:\; \left\{
    \begin{array}{rcl}
    \dot{w}(t)  &  =  &  S\,w(t),  \quad   w_0=w(0) \\
         r(t)       & = &    H  w(t)   \end{array}  \right.   \label{exosyseq}
\ee
 where  $w \in \real^m$ is the state of   exosystem,  and $(x_0,w_0)$  is  an  arbitrary  known  initial condition  for the  nonlinear system \eqref{nonlineq}   with  exosystem \eqref{exosyseq}.

\subsection{Normal  forms for feedback linearizable MIMO systems}

Our methods will assume the system \eqref{nonlineq} is feedback linearizable by state feedback,  so we briefly review the
  standard assumptions needed  to ensure the existence of   suitable linearizing state feedback law  exists \cite{Shastry-99}:
\begin{assumption}  \label{Ass1}
The  origin  is an  equilibrium point of \eqref{nonlineq}.
\end{assumption}
\begin{assumption}  \label{Ass2}
The system \eqref{nonlineq} has a well-defined relative degree vector $(\gamma_1, \dots, \gamma_p)$,  i.e.  there exists $x_0 \in \R^n$ such that
\bn
\item
\bq
 L_{g_i} L_f^k h_i(x)   &   \equiv  &    0, \quad k \in \{0,1,\dots,  \gamma_i -2\}  \\
  L_gL_f^{\gamma-1}  h(x)   &  \neq   &    0
  \eq
  where $g_i$ and  $h_i$ are the  $i$-th components of the vector fields $g$ and  $h$,  for  $i \in \{1,\dots, p\}$.
  \item The matrix
  \be
  A(x) = \left[ \begin{array}{ccc}
               L_{g_1} L_f^{\gamma_1 -1} h_i &  \dots &     L_{g_p} L_f^{\gamma_1 -1} h_i  \\
               \vdots                                 &    \ddots      & \vdots   \\
                 L_{g_1} L_f^{\gamma_p -1} h_i & \dots &   L_{g_p} L_f^{\gamma_p -1} h_i
                 \end{array}\right]
                           \ee
                           is nonsingular at $x= x_0$.
                           \en
                       \end{assumption}
   Under these assumptions,  there  is  neighbourhood $U$ of  $x_0$ on which there exists a change of coordinates
  \be
  \left[ \begin{array}{c} T_1(x)  \\
  \hline
  T_2(x) \end{array}\right]
  = \left[ \begin{array}{c} \eta_1(x)  \\
   \vdots \\
   \eta_{n-\gamma}(x) \\
  \hline
 \xi^1    \\
  \vdots \\
 \xi^p
  \end{array}\right]
  = \left[ \begin{array}{c} \eta  \\
  \hline
  \xi  \end{array}\right]  \label{diffeo}
  \ee
  where  $\gamma = \gamma_1 + \dots + \gamma_p$,   $\frac{\day \eta_i}{\day x}g(x) = 0$ for  $i \in \{1,\dots,n- \gamma\}$,  and for each $j \in \{1,\dots, p\}$,
  \be \xi^j = (\xi_1^j, \xi_2^j, \dots, \xi_{\gamma_j}^j)^T  = (h_j(x), L_f h_j(x), \dots, L_f^{\gamma_j -1} h_j(x))^T
  \ee
  Applying the feedback linearizing control law
    \be u = -A^{-1}(x) \left[\begin{array}{c}
                    L_f^{\gamma_1} h_1   \\
                     \vdots \\
                     L_f^{\gamma_p } h_p \end{array}\right] +   A^{-1}(x) v    \label{mimoulaw}
                     \ee
                     to \eqref{nonlineq} yields the  linear closed-loop system in chain-of-integrator normal  form  coordinates
   \be
   \Sigma_{normal}:\; \left\{ \begin{array}{rcl}
        \dot \eta   & = &  f_o(\eta, \xi)  \\
        \dot \xi & = &  A_c \xi + B_c  v, \quad \xi_0 = \xi(0) \label{normalsys}   \\
           y         &   =  &       C_c \xi
           \end{array} \right.
    \ee
    where $\xi \in \R^\gamma$, $\eta \in \R^{n-\gamma}$,  and $(A_c,B_c,C_c)$ is  a  decoupled MIMO  chain-of-integrator system
   in block diagonal  form  with  \be
 A_c = \blkdiag(A_1,\dots, A_{p}),  \
B_c  = \blkdiag(B_1,\dots, B_{p}), \
 C_c = \blkdiag(C_1,\dots, C_{p}),  \label{Brunsubsys}
 \ee   where,  for each $j \in \{1,\dots,p\}$   each  system is  $(A_j,B_j,C_j)$ is a SISO  chain-of-integrator system  of order $\gamma_j$.
      Finally we also require
        \begin{assumption}  \label{Ass3}
    The zero  dynamics $ \dot \eta = f_0(\eta,0)$ is stable.
    \end{assumption}

    \subsection{Nonovershooting output regulation for Linear systems}
   For the case where  \eqref{nonlineq} is   an  LTI system,  we  have
   \bea
\Sigma_{lin}:\; \left\{ \begin{array}{lclc}
\dot{x}(t)   &  = & A\,x(t)  + B\,u(t)  &  \;\; \; x(0)=x_0 \\
y(t)            &  =  &  C\,x(t)       & \\
\dot{w}(t)  &  =  &  S\,w(t), & \;\; \;  w(0)=w_0 \\
r(t)           & = &    H  w(t)  &
\end{array} \right. \label{Sigmae2}
\eea
For  linear systems,  the problem of  {\it output  regulation by  linear state feedback}  \cite{Saberi-SS-00} can be solved by finding a   control input of  the form
\bea
\label{claw}
u(t)=F\,x(t)+G \,w(t),
\eea
Here $F$ can be any matrix such that $A + B\,F$ is Hurwitz-stable. The  matrices $\Gamma$ and $\Pi$ are obtained  by solving the Sylvester equations
 \bea
\Pi\,S  &  =  &  A\,\Pi+B\,\Gamma \label{Pi1}\\
0  &  =  &  C\,\Pi +H \label{CPi2}
\eea
and finally  $G = \Gamma-F\,\Pi$. Then $u$ as in (\ref{claw}) achieves  output feedback regulation for $\Sigma_{lin}$. Additionally, in  (\cite{SN14}, Theorem 3.1),  it was shown  that the control input \eqref{claw}  yields  a nonovershooting  tracking  response from $(x_0,w_0)$ provided the nominal  system   $\Sigma_{nom} $,  defined  by
 \bea
\Sigma_{nom}:\; \left\{ \begin{array}{lclc}
\dot{\tilde{x}}(t)  &  = &   A\,\tilde{x}(t)     +     B\,\tilde{u}(t),  &   \tilde x(0)= \tilde x_0 \\
\tilde e(t)             &  = &   C\,\tilde{x}(t)  &
\end{array} \right. \label{Sigmanom}
\eea
and  subject to  control  input $\tilde u = F \tilde x$, has a nonovershooting response from initial  condition $\tilde x_0 = x_0 - \Pi w_0$.
\cite{SN10} gave a  range of results for obtaining a nonovershooting  step response from  an LTI system, and in the  next section we adapt these to   linear chain-of-integrator systems.

\section{Nonovershooting natural  response for linear  systems in  chain-of-integrator normal  form} \label{Brunsec}

We consider a  $n$-th order LTI SISO system
\be
\left\{ \begin{array}{lcr}
 \dot x(t) = A\,x(t)+B\,u(t),\;\; x(0)=x_0 \in \real^n, \hfill\cr
 y(t)=C\,x(t),\hfill \end{array} \right.
 \label{sys}
 \ee
 whose  input-output  map is a  chain of  integrators. Thus $A \in \real^{n \times n}$,
 $B \in  \real^{n \times 1}$,  and $C \in   \real^{1 \times n}$,  with
 \be
 A = \left[ \begin{array}{ccccc}
                  0       & 1        & 0         & \dots  & 0 \\
                  0       & 0        & 1             & \ddots  & 0 \\
                 \vdots   & \vdots   & \ddots  & \ddots         & 0  \\
                  0 & 0 & \dots   & 0               &  1           \\
                  0 & 0 & \dots   & 0              & 0
                 \end{array}  \right], \quad
                 B =  \left[ \begin{array}{c}
                  0 \\
                  0  \\
                  \vdots   \\
                   0   \\
                  1
                 \end{array}  \right], \quad
                 C = \left[ \begin{array}{ccccc}
                  1  & 0 & 0  & \dots  & 0
                 \end{array}  \right]   \label{ABCsys}
                 \ee
   Our aim is to  obtain a  feedback  matrix $F$ such  that the state feedback
 control law $u = Fx$ will ensure that the system  natural response $y$ from  $x_0$ converges to $ 0$  without  overshoot.

Firstly let ${\mathcal L} = \{ \lambda_1,  \lambda_{2}, \dots, \lambda_n \}   $  be a set  of  desired   real  stable closed-loop poles with $\lambda_1 < \lambda_2 < \dots < \lambda_n <0$, to be assigned   by   a state feedback matrix $F$. For each $i \in \{1,\dots, n\}$, we solve  the  Rosenbrock  equation
\be
               \left[   \begin{array}{cc}
                    A - \lambda_i I  & B  \\
                            C & 0      \end{array}
                \right]
               \left[   \begin{array}{c}
                    v_i  \\
                    w_i      \end{array}
                \right]
                =   \left[   \begin{array}{c}
                               0  \\ 1                   \end{array}
                \right]  \label{roseneq}
                          \ee
                          where $I$ is the $n \times n$ identity matrix and  the zeros  represent zero  matrices of appropriate dimension.
 Solving     \eqref{roseneq}
yields  vectors ${\mathcal V} = \{v_1, v_2,  \dots, v_n \} \subset \real^n$ and eigendirections  ${\mathcal W} = \{w_1, w_2, \dots, w_n \} \subset \real$  given by
\be
v_i = \left(
        \begin{array}{c}
          1 \\
        \lambda_i \\
        \lambda_i^2 \\
        \vdots \\
        \lambda_i^{n-1}
               \end{array}
      \right), \quad
    w_i =  \lambda_i^n.   \label{viwieq}
      \ee
  We let  $V =  [ v_1 \  v_2 \ \dots \ v_n]$  and $W = [w_1  \ w_2 \ \dots \ w_n]$  and,  following  Moore's  method \cite{Moore76},  obtain   the feedback matrix
\be
F = WV^{-1}   \label{Feq}
\ee
that ensures the eigenstructure of  $A + BF$  consist of eigenvalues  ${\mathcal L}$  and eigenvectors  ${\mathcal V}$.
Introducing $\Lambda =  \diag(\lambda_1,  \dots, \lambda_n)$,   and $\alpha = (\alpha_1, \dots, \alpha_n)^T = V^{-1}x_0 $  we see the  natural response of the closed-loop  system
\be
\left\{ \begin{array}{lcr}
 \dot x(t) = (A+BF)x(t)\;  \quad  x(0)=x_0, \hfill\cr
 y(t)=C\,x(t),\hfill \end{array} \right.
 \label{sys}
 \ee
is  given by
\bq
y(t) &= &  Ce^{(A+BF)t}x_0 \nn  \\
       & = & CVe^{\Lambda t}V^{-1}x_0   \nn  \\
       & = &  \sum_{i=1}^n \   \alpha_i  e^{\lambda_i t}   \label{natresp}
       \eq
       since $Cv_i = 1$ for all $i$,  from \eqref{roseneq}.  Our  next  lemma provides a sufficient condition for $y$ to  not  change sign,  for $t \geq 0$.

    \begin{lemma} \label{NOSlemma}
  For a desired set of real  stable poles $\mathcal L$ and  initial   condition $x_0$, let the natural  response $y$ of the closed-loop  system  be given by \eqref{natresp}. For $ k \in \{1, \dots, n-1\}$, we introduce 
  \be
   c_k:= \left\{ \begin{array}{cc}  1  &  \mbox{ if  $a_k a_n <   0$}  \\
                                            0  &  \mbox{ otherwise}  
                                             \end{array} \right.   \label{ckdef}
\ee 
and  
\be
p(x_0, \mathcal L):=   |\alpha_{n}| +(1-c_{n-1}) |\alpha_{n-1} |   -  \sum_{k=1}^{n-1}  c_k |\alpha_{k}|   \label{p0def} 
\ee
 Then  $y(t) $ does not change sign for   $t \geq 0$  if  $p(x_0, \mathcal L) >0$. 
                              \end{lemma}
      \proof      Assume  $p(x_0, \mathcal L) >0$,  and also that $\alpha_{n}  >0$ and  $ \alpha_{n-1} <0$. Then $c_{n-1} =1$, and for all $t \geq 0$,  by \eqref{p0def}
           \bq
        -(\alpha_1 e^{\lambda_1 t} + \dots  + \alpha_{n-1} e^{\lambda_{n-1} t})
        & \leq &  e^{\lambda_{n-1} t} \sum_{k=1}^{n-1}  c_k |\alpha_{k}|  \nn \\
          & < &  | \alpha_n| e^{\lambda_{n-1} t}   \nn \\
        &   \leq &  \alpha_n e^{\lambda_n t}
        \eq
        Hence
        \be
        0  <  \alpha_1 e^{\lambda_1 t} + \alpha_2 e^{\lambda_2 t} + \dots  +  \alpha_n e^{\lambda_n t} =y(t)
        \ee
        Thus $y(t) >0$ for all $t \geq  0$,  and  hence does not  change sign.  The proof for  $\alpha_{n} >0$  and $  \alpha_{n-1} >0$ is similar.   Corresponding arguments hold when  $\alpha_{n}  <0$  .    \endproof

   In the Appendix,  we  give some more detailed  discussion on  how to choose suitable $\mathcal L$ for   a  given $x_0$, for systems with   dimension $n= 2$ and  $n =3$.

\section{Main Result}

 Under Assumptions \ref{Ass1}-\ref{Ass3}, our task is to  design a  suitable control input  $v$  in \eqref{normalsys} such  that applying  $ u $ in \eqref{mimoulaw}   to \eqref{nonlineq} stabilizes the closed-loop  dynamics and also ensures the  output  tracks the reference signal of  \eqref{exosyseq} without overshoot. We propose the following  algorithm for this problem:
\begin{algorithm} \label{mainalg}
\bn
\item  For each chain-of-integrators subsystem $(A_j,B_j,C_j)$ of  order $\gamma_j$ in \eqref{Brunsubsys}, obtain matrices $\Gamma_j$ and $\Pi_j$ satisfying
 \bea
\Pi_j\,S  &  =  &  A_j\,\Pi+B_j\,\Gamma_j \label{Sylv1}\\
0          &  =  &  C_j\,\Pi_j +H_j \label{Sylv2}
\eea
where $H_j$ denotes the  $j$-th row of $H$ in \eqref{exosyseq}.
Let $\xi_0 = T_2(x_0)$ be  decomposed as
\be \xi_0= (  \xi_0^1,   \xi_0^2, \dots,\  \xi_0^p)^T
\label{Xi0}
\ee  and compute, for each $ j \in \{1, \dots, p\}$,
\be \tilde{\xi}_0^j = \xi_0^j - \Pi_j w_0   \label{txi0j}
\ee
\item Let $\Sigma_{nom}^j$ be the nominal  system in \eqref{Sigmanom} with respect to  $(A_j,B_j,C_j)$,  with initial  condition     $\tilde{\xi}_0^j $. 
Select  candidate closed-loop poles   $\mathcal  L_j = \{\lambda_1^j, \dots, \lambda_{\gamma_j}^j\}$  and solve \eqref{roseneq} to obtain  eigenvectors $\mathcal V_j  = \{v_1^j, v_2^j,  \dots, v_n^j \}$. Compute  $\alpha_j = V_j^{-1}\tilde{\xi}_0^j  $  and  test  $p( \tilde{\xi}_0^j , \mathcal L_j) >0$.  If the test fails,  select  alternative poles. 
   \item Use \eqref{Feq} to obtain the   feedback matrix  $F_j$,  and  compute $G_j = \Gamma_j-F_j\,\Pi_j$.
\item Combine $F = [ F_1^T \  F_2^T \ \dots \ F_p^T ]^T$   and  $G = [ G_1^T \  G_2^T \ \dots \ G_p^T ]^T$  to obtain the control law
\be v = F\xi + Gw    \label{vlaw}
\ee
and  include $v$ in \eqref{mimoulaw} to  obtain  the  feedback linearizing controller  $u$.
\en
\end{algorithm}
Our  main  theorem  sums up  the controller  design method:
\begin{theorem}
Assume the nonlinear system \eqref{nonlineq} meets  Assumptions \ref{Ass1}-\ref{Ass3}.  Let $r$ denote a  reference signal  obtained as the output of the linear exosystem \eqref{exosyseq},  and  let $(x_0,w_0)$ be a known initial condition for \eqref{nonlineq}-\eqref{exosyseq}.  Let   $(\gamma_1, \dots, \gamma_p)$  be  the  relative degree vector of \eqref{nonlineq}, and  let \eqref{normalsys} be its chain-of-integrator normal form under the coordinate change \eqref{diffeo},  with initial condition $\xi_0=  T_2(x_0)$.
Assume that for  each SISO  subsystem $(A_j,B_j,C_j)$ of  \eqref{normalsys},   where $j \in \{1,\dots,p\}$,   the  Sylvester  matrix  equations  \eqref{Sylv1}-\eqref{Sylv2}  admit solutions $\Gamma_j$ and  $\Pi_j$.  Further assume  that for  each subsystem with initial  condition $\tilde \xi_0^j$   given by \eqref{txi0j},  there exist suitable  sets of closed-loop poles   $\mathcal L_j$ such that   $p( \tilde{\xi}_0^j , \mathcal L_j) >0$.  Finally,  let $F$ and $G$ be constructed according to Steps  3 and  4 of Algorithm \ref{mainalg}. Then the feedback linearizing control law \eqref{mimoulaw}, with $v$ given  by \eqref{vlaw}, yields nonovershooting output regulation for the system  \eqref{nonlineq}-\eqref{exosyseq},  from initial condition   $(x_0,w_0)$.
\end{theorem}
{\bf Proof} The decoupled structure of  \eqref{normalsys} ensures that  its $j$-th output component is the output of the $j$-th SISO  subsystem $(A_j,B_j,C_j)$,  which have  been designed to be  nonovershooting. Since the  outputs of  \eqref{nonlineq}-\eqref{exosyseq} under the  control law  $u$ exactly  match those of \eqref{normalsys} under the control  $v$,   we see that  nonovershooting  output regulation is achieved.

\section{Examples}
{\bf Example 4.1}
We   consider the   system    \eqref{nonlineq}  with
\be
f(x) = \left[ \begin{array}{c} x_2 + x_1^2 \\
                                               x_3  \\
                                               x_4  \\
                                               0
                                               \end{array} \right], \
g(x) = \left[ \begin{array}{c}  0 \\
                                                0   \\
                                                 0  \\
                                               1
                                               \end{array} \right], \
    h(x) = x_1, \
     x_0 = \left( \begin{array}{c}1   \\  2 \\  -5 \\   -4\end{array}\right)
            \label{ex1sys}
        \ee
      that was considered in   \cite{Krstic-06}  Example 3, and  \cite{Zhu-13}, Example 2.     The reference signal  was  $r(t) = \cos(t)$,  which can be generated by the exosystem    \eqref{exosyseq} with
      \be
      S=  \bmat{cc}  0 &  1 \\ -1  & 0 \emat, \ H = \bmat{cc} 1 &  0\emat, \  \omega_0 = (1, 0)^T
      \ee
      The coordinate transformation
      \be
      \xi = T(x) = \left( \begin{array}{c}
                   x_1  \\
                   x_2 + x_1^2 \\
                   x_3 +2x_1(x_2+x_1^2) \\
                   x_4 + 2x_1x_3 +(2x_2+6x_1^2)(x_2+x_1^2)
                   \end{array}\right)
                   \ee
                   is   globally singular,  the relative degree  is   4, and  from                 \eqref{mimoulaw},  the  linearizing control input      \be
                   u(x) = -(24x_1^5 +40x_2x_3^3  + 10x_3x_1^2  +(2x_4+16x_2^2)x_1 + 6x_2x_3) + v \label{ex1ulaw}
                   \ee
                   renders the system into  normal  form coordinates  \eqref{normalsys} with $(A_c,B_c,C_c)$ a  fourth-order chain of integrator system.  We apply Algorithm  \ref{mainalg}  as follows.  In Step  1, we solve  the  Sylvester equations   \eqref{Pi1}-\eqref{CPi2} yields
                   \be
                   \Pi = \bmat{cc}
     1  &   0  \\
     0   &  1  \\
      -1  &      0  \\
     0   &  -1
     \emat, \ \Gamma = \bmat{cc}     1   &    0  \emat
     \ee

      We design three feedback matrices  for  a  nonovershooting  system  response from $x_0$ as  follows

     Assuming  the non-zero initial  condition $x_0 = (1, 2, -5, -4)^T$,  we  have
      $\xi_0 = T(x_0) =   (0,  2, -5, 4)^T$  and   $\tilde{\xi}_0  =  \xi_0 - \Pi w_0 = (-1, 2, -4, 4)^T$.   In Step 2,  we select  closed-loop  eigenvalues $\mathcal L$ with $\lambda_1 \in [-6.0 \   -4.5] $, $ \lambda_2 \in [-4.5 \ -3.0]$, $\lambda_3 \in [-3.0 \ -1.5]$ and $\lambda_4 \in [-1.5 \ 0 ]$,  obtain $V$ from \eqref{viwieq} and compute  $\alpha = V^{-1}\tilde \xi_0$.  We use  Lemma  \ref{NOSlemma} to test the suitability of   a given choice of $\mathcal L$. The  choice
      $\mathcal L_1 =  \{          -4.847, \   -4.017, \    -2.432, \    -0.1032\}$ yields  $\alpha =
    ( 0.2468,  -0.3236,   -0.7734,    -0.1499)$,   satisfying $p(x_0, \mathcal L) >0$.    Moving to Step  3, we compute the
      feedback   and  feedforward matrices   $F_1 = -[ 4.89\ 51.6\  42.2 \ 11.4 ]$  and  $G_1 = [  -36.3 \    40.2 ]  $.
    In Step 4  we  form  the control  input $v_1$ in  \eqref{vlaw},  and include  it in   \eqref{ex1ulaw}  to obtain  the control input $u_1$  for \eqref{ex1sys}.

      To   obtain  faster convergence we  also sought  closed-loop  eigenvalues $\mathcal L_2$ with $\lambda_1 \in [-12 \   -9] $, $ \lambda_2 \in [-9 \ -6 ]$, $\lambda_3 \in [-6 \ -3]$ and $\lambda_4 \in [-3  \ 0 ]$,  and $\mathcal L_3$ with $\lambda_1 \in [-16 \   -12] $, $ \lambda_2 \in [-12 \ -8 ]$, $\lambda_3 \in [-8 \ -4]$ and $\lambda_4 \in [-4  \ 0 ]$.  The selections    $\mathcal L_2 =  \{     -10.91, \ -6.55, \   -3.61, \   -2.73 \} $ and
      $\mathcal L_3 =  \{         -15.79, \  -10.20, \  -4.63, \  -3.67\}$  were found to yield  $\alpha$  satisfying   Lemma \ref{NOSlemma},  leading to control matrix pairs  $F_2 =   -[704 \ -625  -192   -23.8  ] $  and    $G_2 = [ 513 \  601]  $;  and  finally $F_3 =  -[2740\   1780 \  394  \ 34]$  and   $G_3 = [ 2347 \   1746]$.

\begin{figure}
\begin{tabular}{cc}
	\includegraphics[height=5.5cm]{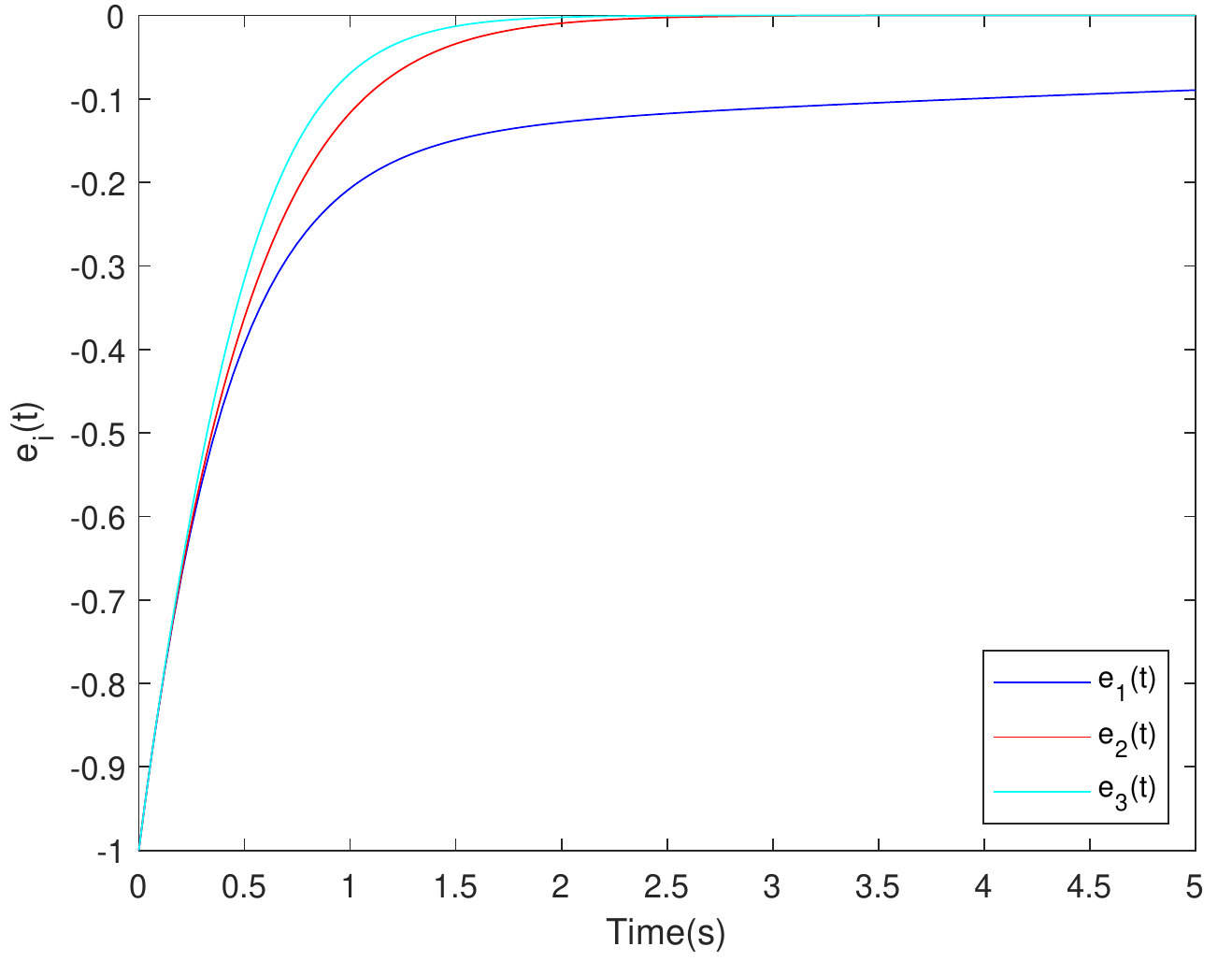} & \ 
	\includegraphics[height=5.5cm]{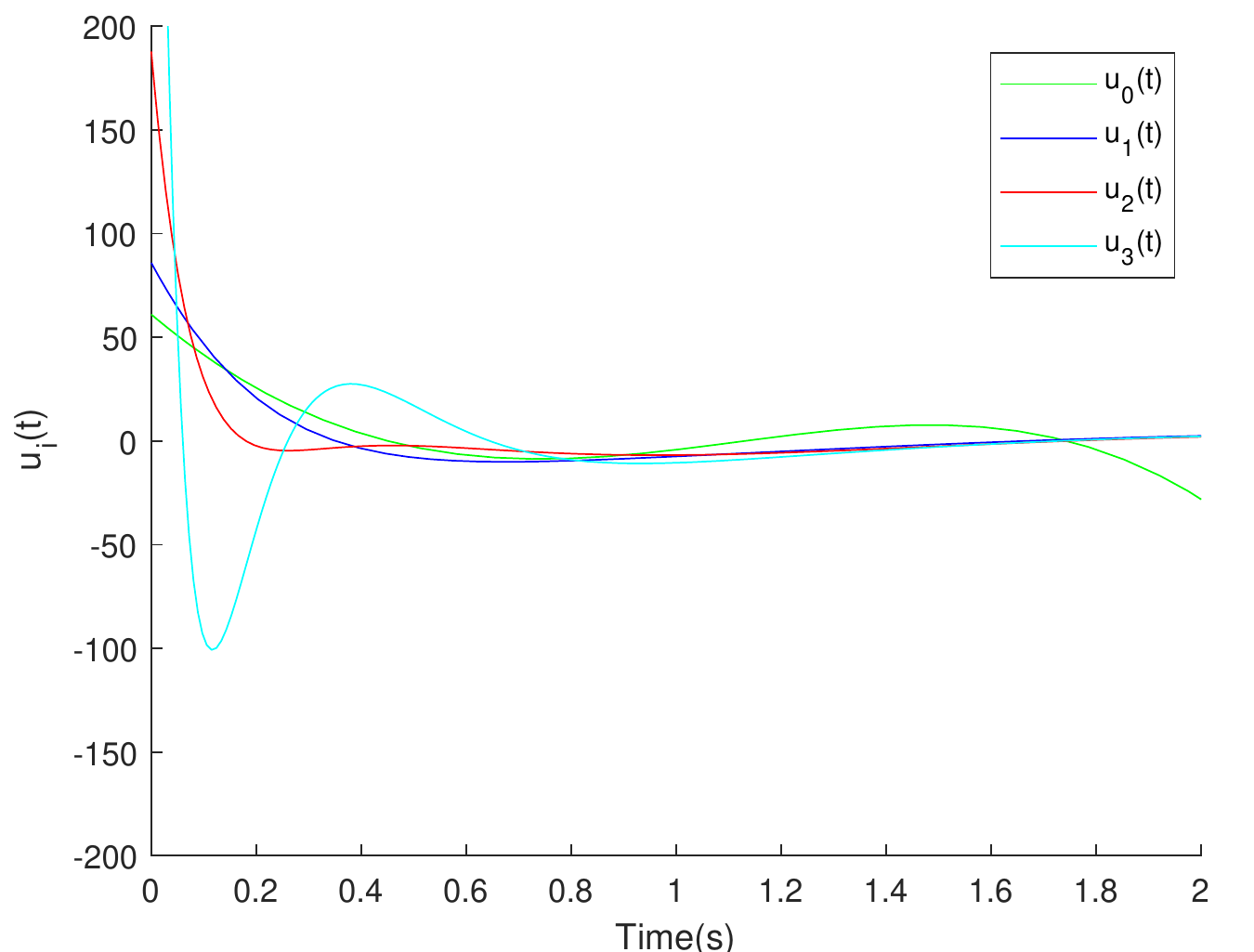}  \\
	 (a)  &  (b) 
	 \end{tabular}
	 	\caption{(a)  Tracking errors  \  (b) Control inputs }
	\label{Fig1}
	
\end{figure}

      Figure  \ref{Fig1}(a)  shows the   error signals $e_i = y_i - r$ arising from each of the  three control inputs;  we  observe the  error signals  do not  change sign and hence in  each   case we achieve a nonovershooting  tracking response.  The  control input signals $u_i$ are shown in Figure \ref{Fig1}(b), including  the  control $u_0$ obtained  from  setting $v =0$. This  input  represents the control  effort required to  linearize the  plant dynamics.

Compared  with  \cite{Krstic-06},  the  advantages of the proposed method lies in its  greater generality,  being applicable to any feedback linearisable  nonlinear system.  Comparing with  \cite{Zhu-13}, both methods assume the system is  feedback linearizable into a chain-of-integrators normal  form. 
 The  method presented here involves the selection of candidate  closed-loop eigenvalues  from  within  specified subintervals of the  negative real   axis, and then testing them for suitability via  Lemma \ref{NOSlemma}.  Thus the  method offers the designer a direct choice between convergence speed and the  magnitude of the controller gains.  Also the selection of the controller gains is apparently rather  more straightforward than the  methods proposed in \cite{Zhu-13},  as it does not involve any further coordinate transformation of the chain-of-integrators system.  The author's own MATLAB implementation of Algorithm \ref{mainalg} on the  chain-of-integrator system requires a modest 40 lines of code, and the search for the feedback  matrix $F_1$  required 20 milliseconds of computation time.

\section*{Appendix}

Here  we consider some  refinements to the results of  Section \ref{Brunsec} that are  available when the chain of integrators system \eqref{ABCsys} are of   second or third order.  Firstly, assume  the state matrices  $(A_c,B_c,C_c) $   are of second order, and   let ${\mathcal L} = \{ \lambda_1,  \lambda_{2}\}   $    be a set  of candidate  stable closed-loop poles with $\lambda_1 < \lambda_2 <0$.  Solving the  Rosenbrock  equation \eqref{roseneq} yields matrices $V$, $W$ and  $F$  in \eqref{Feq}. Letting $x_0 = (x_{01}, x_{02})^T$ be any initial condition, we obtain scalars $\alpha_1,  \ \alpha_2$ such that
    \be
       \left[
        \begin{array}{c}
          \alpha_1   \\
          \alpha_2   \\
               \end{array}
      \right]  =   V^{-1} x_0
      =     \frac{1}{D} \left[
        \begin{array}{c}
          \lambda_2 x_{01} - x_{02}  \\
         - \lambda_1 x_{01} + x_{02}  \\
               \end{array}
      \right]   \label{al1al2}
      \ee
      where $D = \lambda_2 - \lambda_1 >0$. 
               We  introduce
               \be
               q(x_0,\mathcal L)   :=     \alpha_1 \alpha_2 + \alpha_2^2 =    \frac{x_{01}x_{02} - \lambda_1x_{01}^2  }{D} 
           \ee
    By Lemma A.1 of \cite{SN10},  the system has  a nonovershooting response from $x_0$ if $q(x_0,\mathcal L)  > 0$.                  Thus if $x_0$ lies in the first  or third  quadrants of $\real^2$,  we may chose any stable closed-loop poles.
                 If $x_0$ lies in the second  or fourth  quadrants,  we may  choose any $\lambda_2  <0$, but we must choose $\lambda_1 < \frac{x_{02}}{x_{01}}$.

                 Next assume the  state  matrices $(A_c,B_c,C_c) $   are of third order,   let ${\mathcal L} = \{ \lambda_1,  \lambda_{2},\lambda_3\}   $    be a set  of  candidate  stable closed-loop poles with $\lambda_1 < \lambda_2 < \lambda_3 <0$,  and   again obtain  $V$, $W$ and  $F$  in \eqref{Feq}. Letting  $x_0 = (x_{01}, x_{02},x_{03})^T$ be any initial condition, we obtain scalars $\alpha_1,  \ \alpha_2$, $\alpha_3$  such that
       \be
      \left(\begin{array}{c}  \alpha_1 \\ \alpha_2 \\ \alpha_3 \end{array} \right)
       =       V^{-1} x_0
         =  \left(\begin{array}{c} 
  \dsty      \frac{x_{03}\lambda_{2}\,x_{02}\lambda_{3}\,x_{02}+\lambda_{2}\,\lambda_{3}\,x_{01}}{\lambda_{1}^2\lambda_{2}\lambda_{3}}\\
  \dsty   -\frac{x_{03}\lambda_{1}\,x_{02}\lambda_{3}\,x_{02}+\lambda_{1}\,\lambda_{3}\,x_{01}}{\lambda_{1}\lambda_{2}^2\lambda_{3}}\\ 
  \dsty   \frac{x_{03}\lambda_{1}\,x_{02}\lambda_{2}\,x_{02}+\lambda_{1}\,\lambda_{2}\,x_{01}}{\lambda_{1}\lambda_{2}\lambda_{3}^2} \end{array}\right)
      \ee
            From Lemma  \ref{NOSlemma},  \eqref{ckdef} and \eqref{p0def}
            \bq
            c_1            & = &   \frac{f_1(x_0,\mathcal L)f_2(x_0,\mathcal L)  
             }{\left(\lambda_{1}-\lambda_{2}\right)\,{\left(\lambda_{1}-\lambda_{3}\right)}^2\,\left(\lambda_{2}-\lambda_{3}\right)}   \\                     
            c_2  & = &  
                   -\frac{f_1(x_0,\mathcal L)f_3(x_0,\mathcal L)              
             }{\left(\lambda_{1}-\lambda_{2}\right)\,\left(\lambda_{1}-\lambda_{3}\right)\,{\left(\lambda_{2}-\lambda_{3}\right)}^2}  \\
            p(x_0,\mathcal L) 
  & = &      \frac{\left|f_1(x_0,\mathcal L)\right|}{\left(\lambda_{1}-\lambda_{3}\right)\,\left(\lambda_{2}-\lambda_{3}\right)}-\frac{\left|f_3(x_0,\mathcal L)\right|\,\left(2\,c_{2}-1\right)}{\left(\lambda_{1}-\lambda_{2}\right(\,\left(\lambda_{2}-\lambda_{3}\right)}  -\frac{c_{1}\,\left|f_2(x_0,
       \mathcal L)\right|}{\left(\lambda_{1}-\lambda_{2}\right)\,\left(\lambda_{1}-\lambda_{3}\right)}
         \eq  
            
 where                  \bq
                      f_1(x_0,\mathcal L)  & = &
                      x_{03}-(\lambda_{1}+\lambda_{2})\,x_{02}+\lambda_{1}\,\lambda_{2}\,x_{01}  \\
                          f_2(x_0,\mathcal L)  & = &
                          x_{03}-(\lambda_{2}+\lambda_{3})\,x_{02}+\lambda_{2}\,\lambda_{3}\,x_{01}  \\
                            f_3(x_0,\mathcal L)  & = &
                          x_{03}-(\lambda_{1}+\lambda_{3})\,x_{02}+\lambda_{1}\,\lambda_{3}\,x_{01}
                \eq
  Thus  to  ensure a nonovershooting response  from a given $x_0$,  we  need to find suitable $\mathcal L$ to satisfy  the  nonlinear  inequality     $p(x_0,\mathcal L)   >0 $. 

\bibliographystyle{plainnat}

\end{document}